\documentclass[12pt]{article}
\usepackage{amsmath,amssymb}

\newcommand{\ncm}{\newcommand}
\ncm{\aut}{auto\-mor\-phi\-sm} \ncm{\Inn}{\mbox{\rm Inn}} 
\ncm{\Ap}{\mbox{$\overline{\rm Inn}$}} \ncm{\Ext}{\mbox{\rm Ext}} 
\ncm{\Ex}{\mbox{\rm Ex}} \ncm{\OExt}{\mbox{\rm OrderExt}} 
\ncm{\AI}{\mbox{\rm AInn}} \ncm{\HI}{\mbox{\rm HInn($A$)}} 
\ncm{\Aut}{\mbox{\rm Aut}} \ncm{\Mal}{\mbox{$M_{\alpha}$}} 
\ncm{\Aff}{\mbox{${\rm Aff}$}} \ncm{\id}{\mbox{\rm id}} 
\ncm{\Ker}{\mbox{\rm Ker}} \ncm{\BE}{\begin{eqnarray*}} 
\ncm{\EE}{\end{eqnarray*}} \ncm{\lra}{\mbox{$\longrightarrow$}} 
\ncm{\Hom}{\mbox{\rm Hom}} \ncm{\calU}{{\cal U}} \ncm{\el}{\ell} 
\ncm{\ad}{\mbox{\rm ad}} \ncm{\Alg}{\mbox{\rm Alg}} 
\ncm{\Conv}{\mbox{\rm Conv}} \ncm{\D}{{\cal D}} 

\ncm{\cstar}{$C^{*}$-algebra} \ncm{\cstars}{$C^{*}$-algebras} 
\ncm{\ra}{\mbox{$\rightarrow$}} \ncm{\la}{\mbox{$\leftarrow$}} 
\ncm{\hra}{\hookrightarrow} \ncm{\da}{\mbox{$\downarrow$}} 
\ncm{\se}{\mbox{$\searrow$}} \ncm{\del}{\mbox{$\delta$}} 
\ncm{\supp}{\mbox{\rm supp}} \ncm{\Ad}{\mbox{\rm Ad}} 
\ncm{\CAR}{\mbox{$M_{2^{\infty}}$}} \ncm{\ep}{\mbox{$\epsilon > 
0$}} \ncm{\ol}{\overline} \ncm{\Mninf}{\mbox{$M_{n^{\infty}}$}} 
\ncm{\MR}{M. R\o{}rdam} \ncm{\Range}{\mbox{\rm Range}} 
\ncm{\vo}{}
\ncm{\ch}{}
\ncm{\CMP}{Comm. Math. Phys.} \ncm{\add}{} 
\ncm{\tilsig}{\tilde{\sigma}} \ncm{\dist}{{\rm 
dist}}\ncm{\eps}{\epsilon}  \ncm{\calL}{{\mathcal{L}}} 
 
\ncm{\calH}{{\mathcal{H}}} 
\ncm{\lan}{{\langle}}\ncm{\ran}{{\rangle}}

\newtheorem{theo}{Theorem}[section]

\newtheorem{lem}[theo]{Lemma}

\newtheorem{remark}[theo]{Remark}
\newtheorem{definition}[theo]{Definition}
\newtheorem{example}[theo]{Example}
\newtheorem{property}[theo]{Property}

\newenvironment{rem}{\begin{remark} \rm}{\end{remark}}
\newenvironment{pf}{{\it Proof.}}{\hfill$\square$\vspace{3mm}}

\ncm{\R}{\mbox{\bf R}} \ncm{\Z}{\mbox{\bf Z}} \ncm{\T}{\mbox{\bf 
T}} \ncm{\TT}{\T$^{2}$} \ncm{\N}{\mbox{\bf N}} \ncm{\C}{\mbox{\bf 
C}} 

\title{Homogeneity of the pure state space for the separable nuclear \cstars}
\bigskip     
\author{ Akitaka Kishimoto and Sh\^{o}ichir\^{o} Sakai}
\date{April 2001}
\begin{document}
\maketitle 

\begin{abstract} We prove that the pure state space is 
homogeneous under the action of the group of asymptotically inner
automorphisms  for all the separable simple nuclear \cstars. If 
simplicity is not assumed for the \cstars, the set of pure states 
whose GNS representations are faithful is homogeneous for the 
above action.  

\end{abstract} 

 \ncm{\U}{{\mathcal{U}}} 
 \ncm{\F}{{\mathcal{F}}}
 \ncm{\G}{{\mathcal{G}}}
 \ncm{\K}{{\mathcal{K}}}
 \ncm{\M}{{\mathcal{M}}}
 \ncm{\E}{{\mathcal{E}}} 
 \ncm{\B}{{\mathcal{B}}}
  \ncm{\vR}{{\mathcal{R}}}   
 \ncm{\om}{{\omega}}
 \ncm{\al}{{\alpha}}
 \ncm{\Hil}{{\mathcal{H}}}   
 \ncm{\Lin}{{\mathcal{L}}}

\section{Introduction}
If $A$ is a \cstar, an automorphism $\al$ of $A$ is {\em 
asymptotically inner} if there is a continuous family 
$(u_t)_{t\in[0,\infty)}$ in the group  $\U(A)$ of unitaries in $A$ 
(or $A+\C1$ if $A$ is non-unital) such that 
$\al=\lim_{t\rightarrow\infty}\Ad\,u_t$; we denote by $\AI(A)$ the 
group of asymptotically inner automoprphisms of $A$, which is a 
normal subgroup of the group of approximately inner automorphisms. 
Note that each $\al\in\AI(A)$ leaves each (closed two-sided) ideal 
of $A$ invariant. It is shown, in \cite{Pow,Br,FKK}, for a large 
class of separable nuclear \cstars\ that if $\om_1$ and $\om_2$ 
are pure states of $A$ such that the GNS representations 
associated with $\om_1$ and $\om_2$ have the same kernel, then 
there is an $\al\in\AI(A)$ such that $\om_1=\om_2\al$. We shall 
show in this paper that this is the case for all separable nuclear 
\cstars; in particular the pure state space of a separable simple 
nuclear \cstar\ $A$ is homogeneous under the action of $\AI(A)$. 
We do not know of a single example of a separable \cstar\ which 
does not have this property. See \cite{Ef} for some problems on 
this and see \ref{D} and \ref{E} for remarks on the non-separable 
case. 

Choi and Effros \cite{CE3}  have shown that $A$ is nuclear if and 
only if there is a net of pairs $(\sigma_{\nu},\tau_{\nu})$ of 
completely positive (CP) contractons  such that 
$\lim\tau_{\nu}\sigma_{\nu}(x)=x, \ x\in A$, where 
 $$
 A\stackrel{\sigma_{\nu}}{\longrightarrow} N_{\nu} 
 \stackrel{\tau_{\nu}}{\longrightarrow} A
 $$
and $N_{\nu}$ is a finite-dimensional \cstar. When $A$ is a 
non-unital \cstar, $A$ is nuclear if and only if $A+\C1$ is 
nuclear \cite{CE3}. If $A$ is unital, we may assume that both 
$\sigma_{\nu}$ and $\tau_{\nu}$ are unit-preserving. We refer to 
\cite{CE1,CE2} for some other facts on nuclear \cstars. We also 
quote \cite{lance} for a review on the subject. 

Our proof of the homogeneity is a combination of the techniques 
leading up to the above result from \cite{CE3} and the techniques 
from \cite{FKK}. In section 2 we shall show how the homogeneity 
follows from inductive use of Lemma \ref{A} (or \ref{B}), whose 
conclusion is very similar to the properties already used in 
\cite{FKK}; this part follows closely \cite{FKK} and so the proof 
will be sketchy. In section 3 we shall prove Lemma \ref{A} from 
another technical lemma, Lemma \ref{A1}, which shows some 
amenability of the nuclear \cstars; this is the arguments often 
used for individual examples treated in \cite{FKK} and so the 
proof will be again sketchy. Then we will give a proof of Lemma 
\ref{A1}, which constitutes the main body of this paper and uses 
the results and techniques from \cite{CE3}.  

We will conclude this paper, following \cite{FKK}, by generalizing  
Lemma \ref{A1} and then extend the main result, Theorem \ref{C}, 
to show that $\AI(A)$ acts on the pure state space of $A$ {\em 
strongly transitively}. See Theorem \ref{D2} for details.    
  
\section{Homogeneity} 
 \setcounter{theo}{0}   
We first give a main technical lemma, whose conclusion is a 
slightly weaker version of Property 2.6 in \cite{FKK}. We will 
give a proof in the next section. 

\begin{lem} \label{A}
Let $A$ be a nuclear \cstar. Then for any finite subset $\F$ of 
$A$, any pure state $\om$ of $A$ with $\pi_{\om}(A)\cap 
\K(\Hil_{\om})=(0)$, and $\eps>0$, there exist a finite subset 
$\G$ of $A$ and $\delta>0$ satisfying: If $\varphi$ is a pure 
state of $A$ such that $\varphi\sim \om$, and
 $$
 |\varphi(x)-\om(x)|<\delta,\ \ x\in\G,
 $$
then there is a continuous path $(u_t)_{t\in[0,1]}$ in $\U(A)$ 
such that $u_0=1,\ \varphi=\om\Ad\,u_1$, and 
 $$
 \|\Ad\,u_t(x)-x\|<\eps,\ \ x\in\F,\ t\in[0,1].
 $$
\end{lem}

In the above statement, $\pi_{\om}$ is the GNS representation of 
$A$ associated with the state $\om$; $\Hil_{\om}$ is the Hilbert 
space for this representation; $\K(\Hil_{\om})$ is the \cstar\ of 
compact operators on $\Hil_{\om}$; $\varphi\sim\om$ means that 
$\pi_{\varphi}$ is equivalent to $\pi_{\om}$.  We could also 
impose the extra condition that the length of $(u_t)$ is smaller 
than $\pi+\eps$ for the choice of the path $(u_t)$; see Property 
8.1 in \cite{FKK}.  

The following is an easy consequence:
  
\begin{lem}\label{B} 
Let $A$ be a nuclear \cstar. Then for any finite subset $\F$ of 
$A$, any pure state $\om$ of $A$ with $\pi_{\om}(A)\cap 
\K(\Hil_{\om})=(0)$, and $\eps>0$, there exist a finite subset 
$\G$ of $A$ and $\delta>0$ satisfying: If $\varphi$ is a pure 
state of $A$ such that $\ker\pi_{\varphi}=\ker\pi_{\om}$, and
 $$
 |\varphi(x)-\om(x)|<\delta,\ \ x\in\G,
 $$
then for any finite subset $\F'$ of $A$ and $\eps'>0$ there is a 
continuous path $(u_t)_{t\in[0,1]}$ in $\U(A)$ such that $u_0=1$, 
and
 \BE
 |\varphi(x)-\om\Ad\,u_1(x)|<\eps', \ \ && x\in\F',\\
 \|\Ad\,u_t(x)-x\|<\eps,\ \ &&x\in \F.
 \EE
\end{lem}
\begin{pf}
Given $(\F,\om,\eps)$, choose $(\G,\delta)$ as in the previous 
lemma. Let $\varphi$ be a pure state of $A$ such that 
$\ker\pi_{\varphi}=\ker\pi_{\om}$ and 
 $$
 |\varphi(x)-\om(x)|<\delta/2,\ \ x\in\G.
 $$
Let $\F'$ be a finite subset of $A$ and $\eps'>0$ with 
$\eps'<\delta/2$. We can mimic $\varphi$ as a vector state through 
$\pi_{\om}$; by Kadison's transitivity there is a $v\in\U(A)$ such 
that
 $$
 |\varphi(x)-\om\Ad\,v(x)|<\eps',\ \ x\in\F'\cup\G,
 $$
(see 2.3 of \cite{FKK}). Since $|\om\Ad\,v(x)-\om(x)|<\delta,\ 
x\in\G$, we have, by applying Lemma \ref{A} to the pair $\om$ and 
$\om\Ad\,v$, a continuous path $(u_t)$ in $\U(A)$ such that 
$u_0=1$, and
 \BE
 &&\om\Ad\,v=\om\Ad\,u_1,\\
 &&\|\Ad\,u_t(x)-x\|<\eps,\  \ x\in\F.
 \EE
Since $|\varphi(x)-\om\Ad\,u_1(x)|<\eps',\ x\in\F'$, this 
completes the proof.    
\end{pf} 

We shall now turn to the main result stated in the introduction. 
We denote by $\AI_0(A)$ the set of $\al\in\AI(A)$ which has a 
continuous family $(u_t)_{t\in[0,\infty)}$ in $\U(A)$ with $u_0=1$ 
and $\al=\lim\Ad\,u_t$; $\AI_0(A)$ can be smaller than $\AI(A)$ 
(e.g., $\AI_0(A)$ may not contain $\Inn(A)$; see \cite{ER}). 

\begin{theo}\label{C}
Let $A$ be a separable nuclear \cstar. If $\om_1$ and $\om_2$ are 
pure states of $A$ such that $\ker\pi_{\om_1}=\ker\pi_{\om_2}$, 
then there is an $\al\in\AI_0(A)$ such that $\om_1=\om_2\al$. 
\end{theo}          
\begin{pf}
Once we have Lemma \ref{B}, we can prove this in the same way as 
2.5 of \cite{FKK}. We shall only give an outline here. 

Let $\om_1$ and $\om_2$ be pure states of $A$ such that  
$\ker\pi_{\om_1}=\ker\pi_{\om_2}$.

If $\pi_{\om_1}(A)\cap \K(\Hil_{\om_1})\not=(0)$, then  
$\pi_{\om_1}(A)\supset \K(\Hil_{\om_1})$ and $\pi_{\om_1}$ is 
equivalent to $\pi_{\om_2}$. Then by Kadison's transitivity (see, 
e.g., 1.21.16 of \cite{Sak}), there is a continuous path $(u_t)$ 
in $\U(A)$ such that $u_0=1$ and $\om_1=\om_2\Ad\,u_1$. 
                  
Suppose that  $\pi_{\om_1}(A)\cap \K(\Hil_{\om_1})=(0)$, which 
also implies that $\pi_{\om_2}(A)\cap \K(\Hil_{\om_2})=(0)$. 

Let $(x_n)$ be a dense sequence in $A$. 

Let $\F_1=\{x_1\}$ and $\eps>0$ (or $\eps=1$). Let 
$(\G_1,\delta_1)$ be the $(\G,\delta)$ for $(\F_1,\om_1,\eps/2)$ 
as in Lemma \ref{B} such that $\G_1\supset\F_1$. For this 
$(\G_1,\delta_1)$ we choose a continuous path $(u_{1t})$ in 
$\U(A)$ such that $u_{1,0}=1$ and 
 $$
 |\om_1(x)-\om_2\Ad\,u_{1,1}(x)|<\delta_1,\ \ x\in\G_1.
 $$
Let $\F_2=\{x_i,\Ad\,u_{1,1}^*(x_i)\ | \ i=1,2\}$ and let 
$(\G_2,\delta_2)$ be the $(\G,\delta)$ for 
$(\F_2,\om_2\Ad\,u_{1,1},2^{-2}\eps)$ as in Lemma \ref{B} such 
that $\G_2\supset\G_1\cup\F_2$ and $\delta_2<\delta_1$. By \ref{B} 
there is a continuous path $(u_{2t})$ in $\U(A)$ such that 
$u_{2,0}=1$ and 
 \BE
 \|\Ad\,u_{2t}(x)-x\|<2^{-1}\eps,\ \ &&x\in \F_1,\\
 |\om_2\Ad\,u_{1,1}(x)-\om_1\Ad\,u_{2,1}(x)|<\delta_2,\ \ &&x\in\G_2.
 \EE
Let $\F_3=\{x_i,\Ad\,u_{2,1}^*(x_i)\ | \ i=1,2,3\}$ and let 
$(\G_3,\delta_3)$ be the $(\G,\delta)$ for 
$(\F_3,\om_1\Ad\,u_{2,1},2^{-3}\eps)$ as in \ref{B} such that 
$\G_3\supset\G_2\cup\F_3$ and $\delta_3<\delta_2$. By \ref{B} 
there is a continuous path $(u_{3t})$ in $\U(A)$ such that 
$u_{3,0}=1$ and 
 \BE
 \|\Ad\,u_{3t}(x)-x\|<2^{-2}\eps,\ \ &&x\in\F_2,\\    
 |\om_1\Ad\,u_{2,1}(x)-\om_2\Ad(u_{1,1}u_{3,1})(x)|<\delta_3,\ \ 
 &&x\in \G_3.
 \EE
We shall repeat this process. 

Assume that we have constructed $\F_n,\G_n,\delta_n$, and 
$(u_{n,t})$ inductively. In particular if $n$ is even,
 $$
 \F_n=\{x_i,\Ad(u_{n-1,1}^*u_{n-3,1}^*\cdots u_{1,1}^*)(x_i)\ |\ i=1,2,\ldots,n\}
 $$
and $(G_n,\delta_n)$ is the $(\G,\delta)$ for 
$(\F_n,\om_2\Ad(u_{1,1}u_{3,1}\cdots u_{n-1,1}),2^{-n}\eps)$ as in 
\ref{B} such that $\G_n\supset \G_{n-1}\cup\F_n$ and 
$\delta_n<\delta_{n-1}$. And $(u_{n,t})$ is given by \ref{B} for 
$(\F_{n-1},\om_1\Ad(u_{2,1}\cdots u_{n-2,1}),2^{-n+1}\eps)$ and 
for $\F'=\G_n$ and $\eps'=\delta_n$ and it satisfies 
 $$
 |\om_1\Ad(u_{2,1}u_{4,1}\cdots u_{n,1})(x)
 -\om_2\Ad(u_{1,1}\cdots u_{n-1,1})(x)|<\delta_n,\ \ x\in\G_n.
 $$ 
We define continuous paths $(v_t)$ and $(w_t)$ in $\U(A)$ with 
$t\in[0,\infty)$ by: For $t\in [n,n+1]$
 \BE
 &&v_t=u_{1,1}u_{3,1}\cdots u_{2n-1,1}u_{2n+1,t-n},\\
 &&w_t=u_{2,1}u_{4,1}\cdots u_{2n-2,1}u_{2n+2,t-n}.
 \EE
Then, since $\|\Ad\,u_{nt}(x)-x\|<2^{-n+1}\eps,\ x\in \F_{n-1}$, 
we can show that $\Ad\,v_t$ (resp. $\Ad\,w_t$) converges to an 
automorphism $\al$ (resp. $\beta$) as $t\ra\infty$ and that $ 
\om_1\beta=\om_2\al$. Since $\al,\beta\in \AI_0(A)$ and $\AI_0(A)$ 
is a group, this will complete the proof. See the proofs of 2.5 
and 2.8 of \cite{FKK} for details. 
\end{pf} 

The notion of asymptotical innerness for automorphisms may be 
appropriate only for separable \cstars. Because any 
$\alpha\in\AI(A)$ can be obtained as the limit of a {\em sequence} 
in $\Inn(A)$, not just as the limit of a net there. Hence the 
following remark will not be a surprise; it may only suggest that 
we should take $\Ap(A)$ or something bigger than $\AI(A)$ in place 
of $\AI(A)$, in formulating \ref{C} for non-separable \cstars.  

\begin{rem}\label{D}  
There is a unital simple non-separable nuclear \cstar\ $A$ such 
that the pure states space of $A$ is not homogeneous under the 
action of $\AI(A)$. \end{rem} 

We can construct such an example as follows. Let $A$ be a unital 
simple separable nuclear \cstar\ and $\Lambda$ an uncountable set. 
For each finite subset $F$ of $\Lambda$ we set 
$A_F=\otimes_{i\in\Lambda}A_i$ with $A_i\equiv A$ and take the 
natural inductive limit $A_{\Lambda}$ of the net $(A_F)$. Since 
$A_F$ is nuclear, it follows that $A_{\Lambda}$ is nuclear. 

For each $X\subset \Lambda$ we define $A_X$ to be the 
$C^*$-subalgebra of $A_{\Lambda}$ generated by $A_F$ with finite 
$F\subset X$. Note that for each $x\in A_{\Lambda}$ there is a 
countable $X\subset \Lambda$ such that $x\in A_X$.

Let $(u_n)$ be a sequence in $\U(A_{\Lambda})$ such that 
$\Ad\,u_n$ converges to $\alpha\in\Aut(A_{\Lambda})$ in the 
point-norm topology. Since there is a countable subset $X_n\subset 
\Lambda$ such that $u_n\in A_{X_n}$, $\alpha$ is non-trivial only 
on $A_X$, where $X=\cup_nX_n$ is  countable. Thus any $\alpha\in 
\AI(A_{\Lambda})$ has the above property of {\em countable 
support}. 

For each $i\in\Lambda$ let $\om_i$ and $\varphi_i$ be pure states 
of $A_i=A$ such that $\om_i\neq \varphi_i$ and let 
$\om=\otimes_{i\in\Lambda}\om_i$ and 
$\varphi=\otimes_{i\in\Lambda}\varphi_i$. Then it follows that 
$\om$ and $\varphi$ are pure states of $A_{\Lambda}$ and that 
$\om\neq \varphi\alpha$ for any $\alpha\in\AI(A_{\Lambda})$. Hence 
$A_{\Lambda}$ serves as an example for the above remark.

In this case, however, we have an $\alpha\in\Ap(A_{\Lambda})$ such 
that $\om=\varphi\alpha$ (since this is the case for each pair 
$\om_i,\varphi_i$ from \ref{C}) and it may be the case that the 
pure state space of $A_{\Lambda}$ is homogeneous under the action 
of $\Ap(A_{\Lambda})$. 

\begin{rem}\label{E}
There is a unital simple non-separable non-nuclear \cstar\ $A$ 
such that the pure state space of $A$ is not homogeneous under the 
action of $\Aut(A)$. \end{rem}

There are plenty of such \cstars\ at hand. Let $A$ be a factor of 
type II$_1$ or type III with separable predual $A_*$. Then $A$ is 
a unital simple non-separable non-nuclear \cstar\ (see, e.g., 
\cite{lance} for non-nuclearity). Since $A$ contains a 
$C^*$-subalgebra isomorphic to $C_b(\N)\equiv C(\beta\N)$ and 
$\beta\N$ has cardinality $2^c$, the pure state space of $A$ has 
cardinality (at least) $2^c$, where $c$ denotes the cardinality of 
the continuum. (We owe this argument to J. Anderson.) On the other 
hand any $\alpha\in\Aut(A)$ corresponds to an isometry on the 
predual $A_*$, a separable Banach space.  Thus, since the set of 
bounded operators on a separable Banach space has cardinality $c$, 
$\Aut(A)$ has cardinality (at most) $c$. Hence  the pure state 
space of $A$ cannot be homogeneous under the action of  $\Aut(A)$.  

We note in passing that $\AI(A)=\Inn(A)$ for any factor $A$ (or 
any quotient of a factor), since any convergent  sequence in 
$\Aut(A)$ with the point-norm topology converges in norm 
\cite{Ell}. We also note that $\Ap(A)=\Inn(A)$ for any full factor 
\cite{Con1,Sak1}, since then $\Inn(A)$ is closed in $\Aut(A)$ with 
the topology of point-norm convergence in $A_*$ and so is closed 
in $\Aut(A)$ with the topology of point-norm convergence in $A$.

\section{Proof of Lemma \ref{A}} 
 \setcounter{theo}{0}   
If $A$ is a non-unital \cstar, $A$ is nuclear if and only if the 
\cstar\ $A+\C1$ obtained by adjoining a unit is nuclear. Hence to 
prove Lemma \ref{A} we may suppose that $A$ is unital. In the 
following $\U_0(A)$ denotes the connected component of $1$ in the 
unitary group $\U(A)$ of $A$.

\begin{lem}\label{A1}
Let $A$ be a unital nuclear \cstar. Let $\F$ be a finite subset of 
$\U_0(A)$, $\pi$ an irreducible representation of $A$ on a Hilbert 
space $\Hil$, $E$ a finite-dimensional projection on $\Hil$, and 
$\eps>0$. Then there exist an $n\in\N$ and a finite subset $\G$ of 
$M_{1n}(A)$ such that $xx^*\leq1$ and $\pi(xx^*)E=E$ for $x\in\G$, 
and for any $u\in\F$ there is a bijection $f$ of $\G$ onto $\G$ 
with 
 $$\|ux-f(x)\|<\eps.
 $$
\end{lem}

In the above statement, $M_{1n}(A)$ denotes the $1$ by $n$ 
matrices over $A$; if $u\in A$ and $x=(x_1,x_2,\ldots,x_n)\in 
M_{1n}(A)$, 
 \BE
 &&xx^*=\sum_{i=1}^nx_ix_i^*\in A,\\
 &&ux=(ux_1,ux_2,\ldots,ux_n)\in M_{1n}(A).
 \EE
 
We shall first show that Lemma \ref{A1} implies Lemma \ref{A}. 
          
Let $\F$ be a finite subset of $A$, $\om$ a pure state of $A$ with 
$\pi_{\om}(A)\cap \K(\Hil_{\om})=(0)$, and $\eps>0$. Since 
$\U_0(A)$ linearly spans $A$, we may suppose that $\F$ is a finite 
subset of $\U_0(A)$. For $\pi=\pi_{\om}$ and the projection $E$ 
onto the subspace $\C\Omega_{\om}$, we choose an $n\in\N$ and a 
finite subset $\G$ of $M_{1n}(A)$ as in Lemma \ref{A1}. 

We take the finite subset 
 $$
 \{x_ix_j^*\ |\ x\in\G;\ i,j=1,2,\ldots,n\}
 $$
for the subset $\G$ required in Lemma \ref{A}. We will choose 
$\delta>0$ sufficiently small later. Suppose that we are given a 
unit vector $\eta\in\Hil_{\om}$ satisfying
 $$
 |\lan\pi(x_i^*)\eta,\pi(x_j^*)\eta\ran-\lan\pi(x_i^*)\Omega,\pi(x_j^*)\Omega\ran|<\delta
 $$
for any  $x\in\G$ and $i,j=1,2,\ldots,n$, where 
$\Omega=\Omega_{\om}$. Note that 
 $$
 \sum_{j=1}^n\|\pi(x_j^*)\Omega\|^2=\lan\pi(xx^*)\Omega,\Omega\ran=1,
 $$
which implies that $|\lan\pi(xx^*)\eta,\eta\ran-1|<n\delta$.  Thus 
the two finite sets of vectors $S_{\Omega}=\{\pi(x_i^*)\Omega\ |\ 
i=1,\ldots,n;\ x\in\G\}$ and $S_{\eta}=\{\pi(x_i^*)\eta\ |\ 
i=1,\ldots,n;\ x\in\G\}$ have similar geometric properties in 
$\Hil_{\omega}$ if $\delta$ is sufficiently small.  Hence we are 
in a situation where we can apply 3.3 of \cite{FKK}. 

Let us describe how we proceed from here in a simplified case. 
Suppose that the linear span $\Lin_{\Omega}$ of $S_{\Omega}$ is 
orthogonal to the linear span $\Lin_{\eta}$ of $S_{\eta}$ and that 
the map $\pi(x_i^*)\Omega\mapsto \pi(x_i^*)\eta$ and 
$\pi(x_i^*)\eta\mapsto  \pi(x_i^*)\Omega$ extends to a unitary on 
$\Lin_{\Omega}+\Lin_{\eta}$; in particular we have assumed that 
$\lan\pi(x_i^*)\eta,\pi(x_j^*)\eta\ran 
=\lan\pi(x_i^*)\Omega,\pi(x_j^*)\Omega\ran$ for all $i,j$. Since 
$U$ is a self-adjoint unitary, $F\equiv(1-U)/2$ is a projection 
and satisfies that $e^{i\pi F}=U$ on the finite-dimensional 
subspace $\Lin_{\Omega}+\Lin_{\eta}$. By Kadison's transitivity we 
choose an $h\in A$ such that $0\leq h\leq 1$ and 
$\pi(h)|\Lin_{\Omega}+\Lin_{\eta}=F$. We set 
 $$
 \ol{h}=|\G|^{-1}\sum_{x\in\G}xhx^*,
 $$
where 
 $$
 xhx^*=\sum_{i=1}^nx_ihx_i^*.
 $$
Since
 \BE
 \pi(xhx^*)(\Omega-\eta)&=&\sum\pi(x_i)F\pi(x_i^*)(\Omega-\eta),\\ 
 &=&\sum\pi(x_i)\pi(x_i^*)(\Omega-\eta)\\
 &=&\Omega-\eta
 \EE
and $\pi(xhx^*)(\Omega+\eta)=0$, it follows that
 $$
 \pi(\ol{h})(\Omega-\eta)=\Omega-\eta,\ \ 
 \pi(\ol{h})(\Omega+\eta)=0.
 $$
Hence we have that $e^{i\pi \pi(\ol{h})}$ switches $\Omega$ and 
$\eta$. 

On the other hand for $u\in\F$ there is a bijection $f$ of $\G$ 
onto $\G$ such that $\|ux-f(x)\|<\eps,\ x\in\G$. Since
 $$
 u\ol{h}u^*-\ol{h}=|\G|^{-1}\sum_{x\in\G}\{(ux-f(x))hx^*u^*+f(x)h(x^*u^*-f(x)^*)\},
 $$
it follows that $\|u\ol{h}u^*-\ol{h}\|<2\eps$. Thus the path 
$(e^{it\pi \ol{h}})_{t\in[0,1]}$ almost commutes with $\F$ and  is
what is desired. (Since what is required is 
$\om_{\eta}=\om\Ad\,e^{i\pi\ol{h}}$, we may take the path 
$(e^{it\pi(\ol{h}-1/2)})$, whose length is $\pi/2$.) 

If $\Lin_{\eta}$ is not orthogonal to $\Lin_{\Omega}$, we still 
find a unit vector $\zeta\in\Hil_{\om}$ such that 
  $$
 |\lan\pi(x_i^*)\zeta,\pi(x_j^*)\zeta\ran-\lan\pi(x_i^*)\Omega,\pi(x_j^*)\Omega\ran|<\delta
 $$
and such that $\Lin_{\zeta}$ is orthogonal to both $\Lin_{\Omega}$ 
and $\Lin_{\eta}$. Here we use the assumption that 
$\pi_{\om}(A)\cap\K(\Hil_{\om})=(0)$. Then we combine the path of 
unitaries sending $\eta$ to $\zeta$ and then the path sending 
$\zeta$ to $\Omega$ to obtain the desired path. 

The above arguments can be made rigorous in the general case; see 
\cite{FKK} for details. \hfill $\square$

\medskip

We will now turn to the proof of Lemma \ref{A1}, by first giving a 
series of lemmas. The following is an easy version of 3.4 of 
\cite{BKS}. 
 
\begin{lem}\label{A2}
Let $\pi$ be a non-degenerate representation of a \cstar\ $A$ on a 
Hilbert space $\Hil$, $E$ a finite-dimensional projection on 
$\Hil$, $\F$ a finite subset of $A$, and $\eps>0$. Then there is a 
finite-rank self-adjoint operator $H$ on $\Hil$ such that $E\leq 
H\leq 1$ and 
 $$
 \|[\pi(x),H]\|<\eps,\ \ x\in \F.
 $$
\end{lem} 
\begin{pf}
We define finite-dimensional subspaces $V_k,\ k=1,2,\ldots,$ of 
$\Hil$ as follows: $V_1=E\Hil$ and if $V_k$ is defined then 
$V_{k+1}$ is the linear span of $V_k$ and $xV_k,x^*V_k,\ x\in\F$, 
where we have omitted $\pi$ . Then $(V_k)$ is increasing and 
 $$
 x(V_{k+1}\ominus V_k)\subset V_{k+2}\ominus V_{k-1},\ \ x\in \F,
 $$
with $V_0=0$. Denoting by $E_k$ the projection onto $V_k$ we 
define 
 $$
 H_n=\frac{1}{n}\sum_{k=1}^nE_k.
 $$ 
Then $E\leq H_n\leq E_n$. If $x\in\F$, we have, for $\xi\in 
V_{k+1}\ominus V_{k}$, that 
 $$
 (H_nx-xH_n)\xi=(H_n-\frac{n-k}{n})x\xi \in V_{k+2}\ominus 
 V_{k-1}.
 $$
Hence for $\xi\in\Hil$,
 $$
 (H_nx-xH_n)\xi=\sum_{k=0}^{n+1}(H_nx-xH_n)(E_{k+1}-E_k)\xi
 =\sum_{k=0}^{n+1}(H_n-\frac{n-k}{n})x(E_{k+1}-E_k)\xi,
 $$
and thus, by splitting the above sum into three terms, each of 
which is the sum over $k\mod 3=i$ for $i=0,1,2$, and estimating 
each, we reach 
 $$
 \|(H_nx-xH_n)\xi\|\leq \frac{3}{n}\|x\|\|\xi\|.
 $$
This implies that $\|[H_n,x]\|\leq 3/n$ for $x\in\F$.  
\end{pf}

If $\pi$ is a representation of $A$ on a Hilbert space $\Hil$, we 
denote by $\pi_n$ the representation of $M_n\otimes A=M_n(A)$, the 
$n$ by $n$ matrix algebra over $A$, on the Hilbert space 
$\C^n\otimes \Hil$. If $x_i\in A$, then $x_1\oplus 
x_2\oplus\cdots\oplus x_n$ is naturally a diagonal element of 
$M_n(A)$. 

\begin{lem}\label{A3}
Let $\pi$ be a non-degenerate representation of a unital \cstar\ 
$A$ on a Hilbert space $\Hil$, $E$ a finite-rank projection on 
$\Hil$, $\F$ a finite subset of $\U_0(A)$, and $\eps>0$. Then 
there exists an $n\in\N$ such that each $u\in\F$ has a diagonal 
element $\hat{u}= u_1\oplus u_2\oplus\ldots\oplus u_n$ in 
$\U_0(M_n(A))$ satisfying $u_1=u$, $u_n=1$, and 
 $$\|u_i-u_{i+1}\|<\eps/2,\ \ i=1,2,\ldots,n-1.
 $$
Furthermore there exists a finite-rank projection $F$ on 
$\C^n\otimes\Hil$ such that $F\geq E\oplus0\oplus\cdots\oplus0$ 
and
 $$
 \|[\pi_n(\hat{u}),F]\|<\eps,\ \ u\in\F.
 $$ 
\end{lem} 
\begin{pf}
Since $\U_0(A)$ is path-wise connected, the first part is 
immediate.

Let $\delta>0$, which will be specified sufficiently small later. 
By the previous lemma we choose a finite-rank self-adjoint 
operator $H_1$ on $\Hil$ such that $E\leq H_1\leq 1$ and 
 $$
 \|[H_1,u_i]\|<\delta,\ \ i=1,2,\ u\in\F
 $$                                     
where we have omitted $\pi$. Let $E_1$ be the support projection 
of $H_1$ and let $H_2$ be a finite-rank self-adjoint operator on 
$\Hil$ such that $E_1\leq H_2\leq 1$, and 
 $$
 \|[H_2,u_i]\|<\delta,\ \ i=2,3,\ u\in\F.
 $$
In this way we define $H_3,H_4,\ldots,H_{n-1}$ and set 
$H_n=E_{n-1}$, the support projection of $H_{n-1}$. We define an 
operator $F$ on $\C^n\otimes \Hil$ as a tri-diagonal matrix as 
follows:
  \BE
  &&F_{i,i}=H_{i}-H_{i-1},\ \ i=1,\ldots,n, \\
  &&F_{i,i+1}=F_{i+1,i}=\sqrt{H_i(1-H_i)},\ \ i=1,\ldots,n-1,
  \EE
where $H_0=0$. Noting that  $H_iH_{i-1}=H_{i-1}$ and $H_1\geq E$,
it is easy to check that $F$ is a finite-rank projection and $F$ 
dominates $E\oplus0\oplus\cdots\oplus0$. For $u\in\F$, we have 
that 
 \BE
 (\hat{u}F-F\hat{u})_{i,i}&=&[u_i,H_i]-[u_i,H_{i-1}],\\
 (\hat{u}F-F\hat{u})_{i,i+1}&=& [u_i,\sqrt{H_i(1-H_i)}]+ 
 \sqrt{H_i(1-H_i)}(u_i-u_{i+1}).
 \EE
Thus, since $\|\sqrt{H_i(1-H_i)}\|\leq1/2$, the norm of 
$[\hat{u},F]$ is smaller than 
 $$
 \eps/2+2\delta+2\max_i\|[u_i,\sqrt{H_i(1-H_i)}]\|,
 $$
which can be made smaller than $\eps$ for all $u\in\F$ by choosing 
$\delta$ small. 
\end{pf}

When $E$ is a projection on a Hilbert space $\Hil$, we denote by 
$\B(E\Hil)$ the bounded operators on the subspace $E\Hil$.

\begin{lem}\label{A4}
Let $A$ be a unital nuclear \cstar, $\pi$ an irreducible 
representation of $A$ on a Hilbert space $\Hil$, and $E$ a 
finite-rank projection on $\Hil$. Then the identity map on $A$ can 
be approximated by a net of compositions of CP maps
 $$
 A\stackrel{\sigma_{\nu}=\sigma_{\nu}'\oplus\sigma_{\nu}''}{\longrightarrow}
 N_{\nu}\oplus \B(E_{\nu}\Hil)\stackrel{\tau_{\nu}=\tau_{\nu}'+\tau_{\nu}''}
 {\longrightarrow}A,
 $$
where $N_{\nu}$ is a finite-dimensional \cstar, $(E_{\nu})$ is an 
increasing net of finite-rank projections on $\Hil$ such that 
$E\leq E_{\nu}$ and $\lim E_{\nu}=1$, $\sigma_{\nu}'$ and 
$\sigma_{\nu}''$ are unital CP maps such that  
$\sigma_{\nu}''(x)=E_{\nu}\pi(x)E_{\nu},\ x\in A$, and 
$\tau_{\nu}$ is a unital CP map such that 
 \BE
 \pi\tau_{\nu}'(a)E=0,\ \ && a\in N_{\nu},\\ 
 E\pi\tau_{\nu}''(b)E=EbE,\ \ && b\in \B(E_{\nu}\Hil).
 \EE 
\end{lem}     
\begin{pf}
There is a non-degenerate representation $\rho$ of $A$ such that 
$\rho$ is disjoint from $\pi$ and $\rho\oplus\pi$ is a universal 
representation, i.e., $\rho\oplus\pi$ extends to a faithful 
representation of $A^{**}$. Note that 
$(\rho\oplus\pi)(A^{**})=\rho(A)''\oplus\pi(A)''$. 

If the nuclear \cstar\ $A$ is separable, $A^{**}$ is semidiscrete 
\cite{CE1}, which in turn implies that $\vR=\rho(A)''$ is 
semidiscrete. Hence the identity map on $\vR$ can be approximated, 
in the point-weak$^*$ topology, by a net 
$(\tau_{\nu}'\sigma_{\nu}')$ of CP maps on $\vR$, where 
$\sigma_{\nu}'$ (resp. $\tau_{\nu}'$) is a weak$^*$-continuous 
unital CP map of $\vR$ into a finite-dimensional \cstar\ $N_{\nu}$ 
(resp. of $N_{\nu}$ into $\vR$). By denoting $\sigma_{\nu}'\rho$ 
by $\sigma_{\nu}'$ again, we obtain a net of diagrams
 $$
 A\stackrel{\sigma_{\nu}'}{\lra}N_{\nu}\stackrel{\tau_{\nu}'}{\lra}\vR
 $$                                                                
such that $\tau_{\nu}'\sigma_{\nu}'(x)$ converges to $\rho(x)$ in 
the weak$^*$ topology for any $x\in A$.  

If $A$ is separable or not, we have the characterization of 
nuclearity in terms of CP maps \cite{CE3}; there is a net of 
diagrams of unital CP maps:
 $$
 A\stackrel{\sigma_{\nu}'}{\lra}N_{\nu}\stackrel{\tau_{\nu}'}{\lra}A
 $$
such that $N_{\nu}$ is finite-dimensional and 
$\tau_{\nu}'\sigma_{\nu}'(x)$ converges to $x$ in norm for any 
$x\in A$. By denoting $\rho\tau_{\nu}'$ by $\tau_{\nu}'$ again, we 
obtain a net of diagrams:
 $$
 A\stackrel{\sigma_{\nu}'}{\lra}N_{\nu}\stackrel{\tau_{\nu}'}{\lra}\vR
 $$                                                                
as above; actually $\tau_{\nu}'\sigma_{\nu}'(x)$ converges to 
$\rho(x)$ in norm for any $x\in A$.

Since $\pi(A)''=\B(\Hil)$ is semidiscrete, there is such a net of 
CP maps on $\pi(A)''$ as for $\vR$ as well. But we shall construct 
one in a specific way. 

Let $(E_{\nu})$ be an increasing  net of finite-rank projections 
on $\Hil$ such that $E\leq E_{\nu}$ and $\lim E_{\nu}=1$. We 
define $\sigma_{\nu}'':\B(\Hil)\ra \B(E_{\nu}\Hil)$ by 
$\sigma_{\nu}''(x)=E_{\nu}xE_{\nu}$ and 
$\tau_{\nu}'':\B(E_{\nu}\Hil)\ra \B(\Hil)$ by 
$\tau_{\nu}''(a)=a+\omega(a)(1-E_{\nu})$, where $\omega$ is a 
vector state, defined through a fixed unit vector in $E\Hil$. Then 
it is immediate that $(\sigma_{\nu}'',\tau_{\nu}'')$ has the 
desired properties. By denoting $\sigma_{\nu}''\pi$ by 
$\sigma_{\nu}''$ again, we obtain a net of diagrams:
 $$
 A\stackrel{\sigma_{\nu}''}{\lra}\B(E_{\nu}\Hil)\stackrel{\tau_{\nu}''}{\lra}\pi(A)''
 $$  
such that $\tau_{\nu}''\sigma_{\nu}''(x)$ converges to $\pi(x)$ in 
the weak$^*$ topology for any $x\in A$.

We may suppose that we use the same directed set $\{\nu\}$ for 
both $(\sigma_{\nu}',\tau_{\nu}')$  and  
$(\sigma_{\nu}'',\tau_{\nu}'')$. We set 
$\sigma_{\nu}=\sigma_{\nu}'\oplus\sigma_{\nu}''$, 
$M_{\nu}=N_{\nu}\oplus \B(E_{\nu}\Hil)$, and 
$\tau_{\nu}=\tau_{\nu}'+\tau_{\nu}''$. By identifying $A^{**}$ 
with $\vR\oplus\pi(A)''$, we have that
 $$
 A\stackrel{\sigma_{\nu}}{\longrightarrow}M_{\nu}
 \stackrel{\tau_{\nu}}{\longrightarrow}A^{**}
 $$
approximate the identity map on $A$ (in the point-weak$^*$ 
topology), i.e., $\tau_{\nu}\sigma_{\nu}(x)$ converges to $x$ in 
the weak$^*$ topology for any $x\in A$.   

Following \cite{CE3} we approximate $\tau_{\nu}$ by unital CP maps 
of $M_{\nu}$ into $A$. This is done as follows. If $(e_{ij}^k)$ 
denotes a family of matrix units of $M_{\nu}$, $\tau_{\nu}$ is 
uniquely determined by the positive element 
$\Lambda_{\nu}=(\tau_{\nu}(e_{ij}^k))$ in $M_{\nu}\otimes A^{**}$ 
(2.1 of \cite{CE3}). Since $M_{\nu}\otimes A$ is dense in 
$M_{\nu}\otimes A^{**}$ in the weak$^*$ topology, we can, by 
general theory, approximate $\Lambda_{\nu}$ by positive elements 
in $M_{\nu}\otimes A$, in the weak$^*$ topology, which then 
determine CP maps of $M_{\nu}$ into $A$ (see the proof of 3.1 of 
\cite{CE3}).  In particular we approximate $\tau_{\nu}':N_{\nu}\ra 
A^{**}$ by CP maps $\psi':N_{\nu}\ra A$ satisfying 
 $$
 \pi\psi'(a)E=0,\ \ a\in N_{\nu},
 $$
and $\tau_{\nu}'':\B(E_{\nu}\Hil)\ra A^{**}$ by CP maps 
$\psi'':\B(E_{\nu}\Hil)\ra A$ satisfying 
 $$ 
 E\pi\psi''(a)E=EaE,\ \ a\in \B(E_{\nu}\Hil).
 $$
This is indeed possible as shown by using Kadison's transitivity.
Moreover, by taking convex combinations of $\psi'+\psi''$, we may 
assume that $h=\psi'(1)+\psi''(1)$ is close to $1\in A$ in norm. 
By replacing $\psi'$ by $h^{-1/2}\psi'(\,\cdot\,)h^{-1/2}$ etc. we 
may suppose that $\psi=\psi'+\psi''$ is a unital CP map. Since 
$hE=E=Eh$, this does not destroy the above properties imposed on 
$\psi'$ and $\psi''$. 

Restricting $\sigma_{\nu}$ to $A$ and retaining the same symbol 
$\tau$ for the CP maps into $A$ (instead of $\psi$), we now have a 
net of the compositions of unital CP maps: 
 $$
 A\stackrel{\sigma_{\nu}}{\longrightarrow}M_{\nu}
 \stackrel{\tau_{\nu}}{\longrightarrow}A,
 $$                                
which approximates the identity map in the point-weak topology.  

By taking convex combinations of the above CP maps, we will obtain 
such a net which now approximates the identity map in the 
point-norm topology. For example, if $(\lambda_{\nu})$ is such 
that $\lambda_{\nu}\geq0$, $S=\{\nu\ |\ \lambda_{\nu}>0\}$ is 
finite, and $\sum_{\nu}\lambda_{\nu}=1$, then we define 
 $$
 A\stackrel{\phi}{\longrightarrow}(\bigoplus_{\nu\in 
 S}N_{\nu})\oplus \B(E_{\nu_0}\Hil)
 \stackrel{\psi}{\longrightarrow}A,
 $$
where $\nu_0$ is such that $\nu_0\geq \nu,\ \nu\in S$, and
 \BE
 \phi&=&(\oplus_{\nu\in S}\sigma_{\nu}')\oplus \sigma_{\nu_0}'',\\
 \psi&=&(\sum_{\nu\in S}\lambda_{\nu}\tau_{\nu}')
 +(\sum_{\lambda_{\nu}\in S}\lambda_{\nu}\tau_{\nu}''p_{\nu}),
 \EE
with $p_{\nu}:\B(E_{\nu_0}\Hil)\ra \B(E_{\nu}\Hil)$ defined by the 
multiplication of $E_{\nu}$ on both sides. By doing so, the 
properties $\pi\psi'(a)E=0$ and $E\pi\psi''(a)E=EaE$ are still 
retained, where $\psi'$ is the first component of $\psi$ etc. See 
\cite{CE3} for technical details. 
\end{pf}

\begin{lem}\label{A5}
Let $\sigma_{\nu},\tau_{\nu},M_{\nu}=N_{\nu}\oplus\B(E_{\nu}\Hil)$ 
be as in \ref{A4}. For any $\eps>0$ there is a $\delta>0$ such 
that if $u\in \U(A)$ satisfies that 
$\|u-\tau_{\nu}\sigma_{\nu}(u)\|<\delta$, there is a 
$v\in\U(M_{\nu})$ with $\|u-\tau_{\nu}(v)\|<\eps$. 
\end{lem}
\begin{pf}
Suppose that $A$ is represented on a Hilbert space $H$. Since 
$\tau=\tau_{\nu}$ is a unital CP map, by Steinspring's theorem 
there is a representation $\phi$ of $M=M_{\nu}$ on a Hilbert space 
$K$ which contains $H$ such that $\tau(a)=P\phi(a)P,\ a\in M$, 
where $P$ is the projection onto $H$.

If $u\in \U(A)$ satisfies that $\|u-\tau\sigma(u)\|<\delta$, where 
$\sigma=\sigma_{\nu}$ etc., it follows that
 $$
 \tau(\sigma(u)\sigma(u)^*)=P\phi\sigma(u)\phi\sigma(u^*)P
 \geq P\phi\sigma(u)P\phi\sigma(u^*)P\geq (1-2\delta)P.
 $$
Let $b$ denote $\sigma(u)\sigma(u)^*$. Since 
$P\phi(b)(1-P)\phi(b)P=P\phi(b^2)P-(P\phi(b)P)^2\leq 
P-(1-2\delta)^2P$, we have that $\|P\phi(b)(1-P)\|\leq 
2\delta^{1/2}$. Since $[P,\phi(b)]=P\phi(b)(1-P)-(1-P)\phi(b)P$, 
we also have that $\|[P,\phi(b)]\|\leq 2\delta^{1/2}$. For any 
$a\in M$ it follows that $\|\tau(ba)-\tau(b)\tau(a)\|\leq 
2\delta^{1/2}\|a\|$ and  $\|\tau(ba)-\tau(a)\|\leq 
2(\delta^{1/2}+\delta)\|a\|$. 

If $e$ is the spectral projection of $b$ corresponding to 
$[\lambda,1]$ for some $\lambda\in(0,1)$, then $b\leq 
\lambda(1-e)+be$ and
 $$
 (1-2\delta)P\leq P\phi(b)P\leq \lambda P-\lambda P\phi(e)P+P\phi(be)P
 \leq \lambda P-\lambda P\phi(e)P+P\phi(e)P+2(\delta+\delta^{1/2})P.
 $$
Let $\lambda=1-4\delta-2\delta^{1/2}-\delta^{1/4}$. Then the above 
inequality implies that
 $$
 \delta^{1/4}P\leq (4\delta+2\delta^{1/2}+\delta^{1/4})P\phi(e)P,
 $$                                                              
or $\|P-P\phi(e)P\|\leq 4\delta^{3/4}+2\delta^{1/4}$. Hence we 
have that $\|\tau(e)-1\|< 3\delta^{1/4}$ and 
$\|\tau(be)-1\|<3\delta^{1/4}$ for a sufficiently small 
$\delta>0$. Since $be\leq (be)^{1/2}\leq e$, $\tau((be)^{1/2})$ is 
also close to $1$. Since 
$\|\tau(e)-\tau((be)^{1/2})\tau((be)^{-1/2})\|\leq 
\|P\phi((be)^{1/2})(1-P)\| \|(be)^{-1/2}\|<3\delta^{1/8}$, 
$\tau((be)^{-1/2})$ is also close to $1$ (up to the order of 
$\delta^{1/8}$ in this rough estimate); here $(be)^{-1/2}$ is the 
inverse of $(be)^{1/2}$ in $eMe$. 

We now define a unitary $v$ in $M$ by $v=(be)^{-1/2}\sigma(u)+y$, 
where $y$ satisfies that $yy^*=1-e$ and 
$y^*y=1-\sigma(u)^*(be)^{-1}\sigma(u)$. Since 
$(be)^{-1/2}\sigma(u)\sigma(u)^*(be)^{-1/2}=e$, $v$ is indeed a 
unitary. Since $\tau(y)\tau(y^*)\leq \tau(yy^*)=\tau(1-e)\leq 
3\delta^{1/4}$, $\|y\|$ is of the order of $\delta^{1/8}$. Since 
$\tau((be)^{-1/2}\sigma(u))$ is close to 
$\tau((be)^{-1/2})\tau(\sigma(u))$ up to the order of 
$\delta^{1/16}$, we can conclude that 
$\|\tau(v)-\tau(\sigma(u))\|$ is close to zero up to the order of 
$\delta^{1/16}$. 
\end{pf}

When $(X,d)$ is a metric space, $S\subset X$, and $\eps>0$, we 
call $S$ an $\eps$-net if $\cup_{x\in S}B(x,\eps)=X$, where 
$B(x,\eps)=\{y\in X\ |\ d(x,y)<\eps\}$. When $X$ has a finite 
$\eps$-net, we denote by $N(X,\eps)$ the minimum of orders over 
all the finite $\eps$-nets. If $X$ is compact, then $N(X,\eps)$ is 
well-defined for any $\eps>0$. 

\begin{lem}\label{A6}                    
Let $(X,d)$ be a compact metric space. If $S_1$ and $S_2$ are 
$\eps$-nets consisting $N(X,\eps)$ points, then there is a 
bijection $f$ of $S_1$ onto $S_2$ such that $d(x,f(x))<2\eps,\ 
x\in S_1$.  
\end{lem} 
\begin{pf}
Let $\F$ be a non-empty subset of $S_1$ and set
 $$
 \G=\{y\in S_2\ |\ B(y,\eps)\cap \cup_{x\in 
 \F}B(x,\eps)\not=\emptyset\}.
 $$
Since $\cup_{x\in\F}B(x,\eps)\subset \cup_{x\in\G}B(x,\eps)$, it 
follows that $\G\cup S_1\setminus \F$ is an $\eps$-net and that 
the order of $\G$ is greater than or equal to the order of $\F$. 
Then by the matching theorem we can find a bijection $f$ of $S_1$ 
onto $S_2$ such that $f(x)\in\{y\in S_2\ |\ B(x,\eps)\cap 
B(y,\eps)\not=\emptyset\}$. \end{pf}
                             
\noindent {\em Proof of Lemma \ref{A1}} 
 Let $\pi$ be an irreducible representation of the unital nuclear \cstar\ 
$A$ on a Hilbert space $\Hil$, $E$ a finite-rank projection on 
$\Hil$, $\F$ a finite subset of $\U_0(A)$, and $\eps>0$.

We apply Lemma \ref{A3} to this situation. Thus there exist an 
$n\in\N$ and a finite-rank projection $F$ on $\C^n\otimes\Hil$ 
such that 
 \BE
 &&F\geq E\oplus0\oplus\cdots\oplus0,\\ 
 &&\|[F,\pi_n(\hat{u})]\|<\eps,\ \ u\in\F,
 \EE                               
where $\pi_n$ denotes the natural extension of $\pi$ to a 
representation of $M_n\otimes A$ on $\C^n\otimes\Hil$; hereafter 
we shall simply denote $\pi_n$ by $\pi$. Let $F_0$ be a 
finite-rank projection on $\Hil$ such that $F\leq 1\otimes F_0$.    

By Lemma \ref{A4} we find a net of diagrams
 $$
 A\stackrel{\sigma_{\nu}=\sigma_{\nu}'\oplus\sigma_{\nu}''}{\longrightarrow}
 N_{\nu}\oplus \B(E_{\nu}\Hil)
 \stackrel{\tau_{\nu}=\tau_{\nu}'+\tau_{\nu}''}{\longrightarrow}A
 $$
with $F_0$ in place of $E$ as described there; in particular 
$F_0\leq E_{\nu}$. We take tensor product of these diagrams with 
$M_n$; denoting $\id_n\otimes\sigma_{\nu}$ by the same symbol 
$\sigma_{\nu}$ etc., we obtain 
 $$
 M_n\otimes A\stackrel{\sigma_{\nu}=\sigma_{\nu}'\oplus\sigma_{\nu}''}{\longrightarrow}
 M_n\otimes N_{\nu}\oplus M_n\otimes\B(E_{\nu}\Hil)
 \stackrel{\tau_{\nu}=\tau_{\nu}'+\tau_{\nu}''}{\longrightarrow}M_n\otimes 
 A.
 $$
Noting that $F\in M_n\otimes\B(E_{\nu}\Hil)=\B(\C^n\otimes 
E_{\nu}\Hil)$, we denote 
 $$
 V_{\nu}=\U( M_n\otimes N_{\nu}\oplus M_n\otimes\B(E_{\nu}\Hil)\cap \{F\}'), 
 $$                                                                     
which is a compact group. Since $(1\otimes 
F_0)\pi\tau_{\nu}'(v)=0$ and $(1\otimes 
F_0)\pi\tau_{\nu}''(v)(1\otimes F_0)=(1\otimes F_0)v(1\otimes 
F_0)$, we have that  for each $v\in V_{\nu}$ 
 \BE
 F\pi(\tau_{\nu}(v)\tau_{\nu}(v^*))F
 &=&F(1\otimes F_0)\pi(\tau_{\nu}(v)\tau_{\nu}(v^*))(1\otimes F_0)F,\\
 &=&F(1\otimes F_0)\pi(\tau_{\nu}''(v)\tau_{\nu}''(v^*))(1\otimes F_0)F,\\  
 &=& F(1\otimes F_0)v(1\otimes F_0)v^*(1\otimes F_0)F \\
  &&+ F(1\otimes F_0)\pi(\tau_{\nu}''(v))(1\otimes(1-F_0))\pi(\tau_{\nu}''(v^*))(1\otimes F_0)F.
  \EE
Since the first term is $F$ as $[F,v]=0$, the second term must be 
zero. Hence it follows that 
 $$
 F\pi(\tau_{\nu}(v)\tau_{\nu}(v)^*)F=F,
 $$
which implies that
 $$
 \pi(\tau_{\nu}(v)\tau_{\nu}(v)^*)F=F.
 $$                                  
By multiplying $E\oplus0\oplus\cdots\oplus0$ from the right we 
have that
 $$
 \sum_{j,k}\pi(\tau_{\nu}(v_{1j})\tau_{\nu}(v^*_{kj}))F_{k1}E=E.
 $$
Since $F\geq E\oplus0\oplus\cdots\oplus0$, we have that 
$F_{k1}E=0$ for $k\neq1$. Thus it follows that for $v\in V_{\nu}$, 
 $$
 \sum_{j=1}^n\pi(\tau_{\nu}(v_{1j})\tau_{\nu}(v_{1j}^*))E=E.
 $$

By Lemma \ref{A5} (applied to $M_n\otimes A$ instead of $A$) we 
choose $\nu$ such that each $u\in\F$ has a unitary $\hat{u}'\in  
M_n\otimes N_{\nu}\oplus M_n\otimes\B(E_{\nu}\Hil)$ such that
 $$
 \|\tau_{\nu}(\hat{u}')-\hat{u}\|\approx0
 $$
as well as
 $$
 \|\tau_{\nu}\sigma_{\nu}(\hat{u})-\hat{u}\|\approx0.
 $$
Since
 \BE
 &&(1\otimes F_0)\hat{u}'(1\otimes F_0)
 = (1\otimes F_0)\pi(\tau_{\nu}''(\hat{u}'))(1\otimes F_0) \\
 &&\approx(1\otimes F_0)\pi(\tau_{\nu}(\hat{u}'))(1\otimes F_0)
 \approx (1\otimes F_0)\pi(\hat{u})(1\otimes F_0),
 \EE
we have that
 $$
 \pi(\hat{u})F\approx F \pi(\hat{u})F\approx F\hat{u}'F\approx 
 \hat{u}'F.
 $$
By choosing $\nu$ sufficiently large, we may assume that
 $$
 \|[\hat{u}',F]\|<\eps,\ \ u\in \F.
 $$                                
By taking the unitary part of the polar decomposition of 
$w=F\hat{u}'F+(1-F)\hat{u}'(1-F)$, we may assume that
 $$
 [\hat{u}',F]=0,\ \ u\in\F.
 $$
Since $\|w-\hat{u}'\|<2\eps$, we can estimate that
 $$
 \|\tau_{\nu}(\hat{u}')-\hat{u}\|<3\eps,\ \ u\in\F.
 $$  
Since $\|\tau_{\nu}(\hat{u}')\tau_{\nu}(\hat{u}')^*-1\|<6\eps$, we 
have that for any $v\in V_{\nu}$,
 $$
 \|\tau_{\nu}(\hat{u}'v)-\tau_{\nu}(\hat{u}')\tau_{\nu}(v)\|<(12\eps)^{1/2}<4\eps^{1/2}.
 $$ 
(See the proof of \ref{A5}.) Hence for $v\in V_{\nu}$
 $$
 \|\hat{u}\tau_{\nu}(v)-\tau_{\nu}(\hat{u}'v)\|<3\eps+4\eps^{1/2}, 
 \ \ u\in\F.
 $$

We choose an $\eps$-net $\G'$ of $V_{\nu}$ consisting of 
$N(V_{\nu},\eps)$ points and set 
 $$
 \G=\{(\tau_{\nu}(v_{11}),\tau_{\nu}(v_{12}),\ldots,\tau_{\nu}(v_{1n}))\ |\ v\in\G'\}. 
 $$ 
Since $\hat{u}'\G'$ is also an $\eps$-net of $V_{\nu}$ for 
$u\in\F$, Lemma \ref{A6} gives a bijection $f$ of $\G'$ onto $\G'$ 
such that 
 $$
 \|\hat{u}'v-f(v)\|<2\eps,\ \ v\in\G'.
 $$
Hence for each $u\in\F$ there is a bijection $f$ of $\G'$ onto 
$\G'$ such that 
 $$
 \|\hat{u}\tau_{\nu}(v)-\tau_{\nu}(f(v))\|<5\eps+4\eps^{1/2},
 $$                                                           
which implies that regarding $f$ as a map of $\G$ onto  $\G$,
 $$
 \|ux-f(x)\|<5\eps+4\eps^{1/2},\ \ x\in\G.
 $$
This completes the proof. \hfill $\square$ 

\medskip

In Lemma \ref{A4} we could handle a mutually disjoint finite 
family of irreducible representations instead of just one. By 
doing so we can derive:

\begin{lem}\label{D1}
Let $A$ be a unital nuclear \cstar. Let $\F$ be a finite subset of 
$\U_0(A)$, $\pi$ a  representation of $A$ on a Hilbert space 
$\Hil$ such that $\pi=\oplus_{i=1}^k\pi_k$ with $(\pi_i)_{i=1}^k$ 
a mutually disjoint family of irreducible representations of $A$, 
$E$ a finite-dimensional projection on $\Hil$, and $\eps>0$. Then 
there exist an $n\in\N$ and a finite subset $\G$ of $M_{1n}(A)$ 
such that $xx^*\leq1$ and $\pi(xx^*)E=E$ for $x\in\G$, and for any 
$u\in\F$ there is a bijection $f$ of $\G$ onto $\G$ with 
 $$\|ux-f(x)\|<\eps.
 $$
\end{lem} 

A straightforward generalization of \ref{A4} would require that
$E\in\pi(A)''$ in the above statement. But, since any finite-rank 
projection on $\Hil$ is dominated by such a one in $\pi(A)''$, we 
did not need it. 

By having this at hand we can derive a stronger version of Lemma 
\ref{A} and then strengthen Theorem \ref{C}. For example we will 
obtain:

\begin{theo}\label{D2}
Let $A$ be a separable nuclear \cstar. If $(\om_i)_{1\leq i\leq 
n}$ and $(\varphi_i)_{1\leq i\leq n}$ are finite sequences of pure 
states of $A$ such that $(\om_i)$ (resp. $(\varphi_i)$) are 
mutually disjoint and $\ker_{\om_i}=\ker_{\varphi_i}$ for all $i$, 
then there is an $\alpha\in\AI_0(A)$ such that 
$\om_i=\varphi_i\alpha$ for all $i$. \end{theo}

We will have to use a general form of Kadison's transitivity for 
the proofs of the above results as in \cite{Sak}. See Section 7 of 
\cite{FKK} for details and for other consequences.    

We do not know whether we could take an arbitrary non-degenerate 
representation of $A$ for $\pi$ in Lemma \ref{D1} (perhaps by 
weakening the requirement $\pi(xx^*)E=E$ by 
$\|\pi(xx^*)E-E\|<\eps$). If this were the case, we would obtain a 
new characterization of nuclearity which manifests a close 
connection with  amenability of $A$ (cf. \cite{Con2,Haa0,Pat}).

\medskip
\small
\begin{flushright}
 Department of Mathematics, Hokkaido University, Sapporo, Japan 060-0810 \\ 
 5-1-6-205, Odawara, Aoba-ku,  Sendai, Japan 980-0003 
\end{flushright}


\small
\begin{thebibliography}{99}
\bibitem{Br} O. Bratteli, Inductive limits of finite-dimensional 
\cstars, Trans. Amer. Math. Soc. 171 (1972), 195--234. 
  
\bibitem{BKS} N.P. Brown, K. Dykema, and D. Shlyakhtenko, 
Topological entropy of free product automorphisms, preprint.
\bibitem{CE1} M-D. Choi and E.G. Effros, Separable nuclear \cstars\ 
and injectivity, Duke Math. J. 43 (1976), 309--322.
\bibitem{CE2} M-D. Choi and E.G. Effros, Nuclear \cstars\ and 
injectivity: The general case, Indiana Univ. Math. J. 26 (1977), 
443--446.
\bibitem{CE3} M-D. Choi and E.G. Effros, Nuclear \cstars\ and the 
approximation property, Amer. J. Math. 100 (1978), 61--79.   


\bibitem{Con1} A. Connes, Almost periodic states and factors of 
type III$_1$, J. Funct. Anal. 16 (1974), 415--445.

\bibitem{Con2} A. Connes, On the cohomology of operator algebras, J. 
Funct. Anal. 28 (1978), 248--253.

\bibitem{Ef} E.G. Effros, On the structure theory of \cstars: some old and new
problems, in: Proceedings of symposia in pure mathematics 38 
(1982) part 1, edited by R.V. Kadison, pages 19--34.  

\bibitem{Ell} G.A. Elliott, Convergence of automorphisms in certain \cstars, 
J. Funct. Anal. 11 (1972), 204--206. 
\bibitem{ER} G.A. Elliott and M. R\o{}rdam, The automorphism group 
of the irratinal rotation \cstar, Commun. Math. Phys. 155 (1993), 
3--26. 
\bibitem{FKK} H. Futamura, N. Kataoka, and A. Kishimoto, Homogeneity 
of the pure state space for separable \cstars, to appear in 
Internat. J. Math.  
\bibitem{Haa0} U. Haagerup, All nuclear \cstars\ are amenable, Invent. Math. 74 (1983),
305--319.    

 

\bibitem{lance} E.C. Lance, Tensor products and nuclear \cstars, 
in: Proceedings of symposia in pure mathematics 38 (1982) part 1, 
edited by R.V. Kadison, pages 379--399.  
\bibitem{Pat} A.T. Paterson, Nuclear \cstars\ have amenable unitary groups,
Proc. Amer. Math. Soc. 114 (1992), 719--721. 
 \bibitem{Pow} R.T. Powers, Representations of uniformly hyperfinite 
algebras and their associated von Neumann rings, Ann. of Math. 86 
(1967), 138--171. 
\bibitem{Sak1} S. Sakai, On automorphism groups of II$_1$-factors, 
T\^{o}hoku Math. J. 26 (1974), 423--430.
\bibitem{Sak} S. Sakai, {\em $C^*$-algebras and $W^*$-algebras}, 
Classics in Math., Springer, 1998. 
\end{thebibliography}
\end{document}